\documentclass[proc]{edpsmath}
\usepackage[dvipdf]{graphicx}
\usepackage[latin1]{inputenc}
\usepackage{latexsym}
\usepackage{amsmath,amssymb,amsfonts,amsthm}
\usepackage{xcolor}
\usepackage{mathrsfs}
\usepackage[colorlinks=true,citecolor=red,linkcolor=blue,urlcolor=blue,pdfstartview=FitH]{hyperref}

\usepackage{tikz}\usetikzlibrary{arrows,calc,decorations.text,decorations.pathreplacing}
\usepackage{qtree}
\usepackage{times}
\usepackage{color}
\usepackage{array}
\usepackage{bbold}
\usepackage{graphicx}
\usepackage{mathrsfs}

\usepackage{amsmath, latexsym, amsfonts, amssymb, amsthm, amscd}
\usepackage{graphicx,epsf,psfrag}

\usepackage{hyperref}
\hypersetup{colorlinks=true,breaklinks=true,urlcolor=blue,linkcolor=blue,bookmarksopen=true} 

\newcommand{\Exp}{\mathbb{E}} 
\newcommand{\Wass}{\mathrm{W}} 

\newcommand{\dd}{\,\text{\rm d}}             

\newcommand{\ER}{\mathbb {R}}
\newcommand{\R}{\mathbb {R}}

\newcommand{\PE}{\mathbb {P}}
\newcommand{\ES}{\mathbb{E}}
\newcommand{\abs[1]}{\lvert#1\rvert}

\newcommand{\bqn}{\begin{equation}}
\newcommand{\bqne}{\begin{equation*}}
\newcommand{\eqn}{\end{equation}}
\newcommand{\eqne}{\end{equation*}}

\newcommand*{\N}{\mathbb{N}}

\newcommand*{\E}{\mathbb{E}}

\newcommand*{\e}{\text{e}}

\newcommand*{\eps}{\varepsilon}
\newcommand*{\prob}{\mathbb{P}}

\begin{document}

\title{Long time behavior of Markov processes and beyond}
\thanks{The authors thanks the Agence Nationale de la Recherche PIECE 12-JS01-0006-01 
for its financial support.}
%
\author{Florian Bouguet}\address{Institut de Recherche Math\'ematique 
de Rennes (UMR 6625 CNRS), Universit\'e de Rennes 1, Campus de Beaulieu, F-35042 Rennes 
Cedex, France.}

\author{Florent Malrieu}\address{Laboratoire de Math\'ematiques et Physique Th\'eorique 
(UMR CNRS 7350), F\'ed\'eration Denis Poisson (FR CNRS 2964), Universit\'e 
Fran\c cois-Rabelais, Parc de Grandmont, 37200 Tours, France.}

\author{Fabien Panloup}\address{Institut de Math\'ematiques de Toulouse (UMR CNRS 5219),
Universit\'e Paul Sabatier \& INSA Toulouse, 135, av. de Rangueil, F-31077 Toulouse 
Cedex 4, France.}

\author{Christophe Poquet}\address{Dipartimento di Matematica, Universit\`a degli 
studi di Roma Tor Vergata, Via della Ricerca Scientifica I-00133 Roma, Italy.}

\author{Julien Reygner}\address{Laboratoire de Physique (UMR CNRS 5672), 
\'Ecole Normale Sup\'erieure de Lyon, 46 all\'ee d'Italie, Lyon, 69364, France.}

\begin{abstract} This note provides several recent progresses in the study of long 
time behavior of Markov processes. The examples presented below are 
related to other scientific fields as PDE's, physics or biology. The involved 
mathematical tools as propagation of chaos, coupling, functional inequalities, 
provide a good picture of the classical methods that furnish quantitative rates 
of convergence to equilibrium.
\end{abstract}

\begin{resume} Cet article présente plusieurs progrès récents dans l'étude du 
comportement en temps long de certains processus de Markov. Les exemples 
présentés ci-dessous sont motivés par différentes applications issues de la physique 
ou de la biologie. Les outils mathématiques employés, propagation du chaos, 
couplage, inégalités fonctionnelles, couvrent un large spectre des techniques 
disponibles pour obtenir des comportements en temps long quantitatifs. 
\end{resume}

\maketitle

\section*{Introduction}

This note gathers several progresses in the study of the long time behavior of Markovian 
(and non Markovian) processes. The first section is dedicated to the study of stochastic 
differential equation driven by a fraction Brownian motion with a non constant diffusion 
matrix. This process is not Markovian but one can successfully adapt the coupling strategy 
to get quantitative long time estimates. In the second section, a piecewise deterministic Markov 
process arising linked to a stochastic algorithm is studied thanks to clever couplings of the 
paths. The last two sections stress the fruitful links between mean field interacting 
particle systems and non linear parabolic partial differential equations.

\section{Rate of convergence to equilibrium for fractional SDEs}
This part, which is a short version of \cite{fontbona-panloup}, is devoted to the 
problem of the estimation of the rate of convergence to equilibrium of stochastic 
differential equations (SDEs) driven by a fractional Brownian motion (fBm). In this 
highly not Markovian setting, this problem has been first investigated by Hairer~\cite{hairer} 
who proved that, in the additive setting, the rate of convergence in total variation can be 
upper-bounded by $C t^{-\rho_H}$ where $\rho_H$ is a positive number depending on 
the Hurst parameter $H$ of the fBm. 
But to our knowledge, there is no result in the multiplicative setting. Hence, we focus on 
the generalization of the existing additive results to the multiplicative case and  obtain an 
extension of \cite{hairer}  when $H>1/2$.  The main novelty  of this work  is a development 
of some Foster-Lyapunov  techniques in this non-Markovian setting, which  allows us to 
put in place an asymptotic coupling scheme as that of  \cite{hairer} without resorting 
to deterministic contracting properties.

\subsection{Introduction} 

We deal   with an $\ER^d$-valued process $(X_t)_{t\ge0}$ 
which is a solution to the following SDE
\begin{equation}\label{fractionalSDE0}
dX_t=b(X_t)dt+ \sigma(X_t)dB_t
\end{equation}
where $b:\ER^{d}\rightarrow\ER^d$ and $\sigma:\ER^{d}\rightarrow \mathbb{M}_{d,d} $ 
are (at least) continuous functions, and where $ \mathbb{M}_{d,d}$ is the  set of $d\times d$ 
real matrices. In~\eqref{fractionalSDE0}, $(B_t)_{t\ge0}$ is  a $d$-dimensional fBm with 
Hurst parameter $H\in(\frac{1}{2},1)$. Note that under some Hölder regularity assumptions 
on the coefficients (see $e.g.$ \cite{Nualart02} for background), (strong) existence and 
uniqueness hold for the solution to \eqref{fractionalSDE0} starting from $x_0\in\ER^d$. 

Introducing the Mandelbrot-Van Ness representation of the fractional Brownian motion, 
\begin{equation}\label{eq:mandel}
B_t=\alpha_H\int_{-\infty}^{0} (-r)^{H-\frac{1}{2}} \left(dW_{r+t}-dW_r\right),\quad t\ge0,
\end{equation}
where $(W_t)_{t\in\ER}$ is a two-sided $\ER^d$-valued Brownian motion and $\alpha_H$ 
is a normalization coefficient depending on $H$, $(X_t, (B_{s+t})_{s\le 0})_{t\ge0}$ can 
be realized through a Feller transformation $({\cal Q}_t)_{t\ge0}$ on the product space 
$\ER^d\times{\cal W}$ where ${\cal W}$ is a suitable Hölder-type space on $(-\infty,0]$ 
(we refer to \cite{hairer2} for more rigorous background on this topic). In particular,
an initial distribution of this dynamical system is a distribution $\mu_0$ on 
$\ER^d\times{\cal W}$. With probabilistic words, an initial distribution is the distribution 
of  a couple $(X_0, (B_s)_{s\le0})$ where $(B_s)_{s\le0}$ is an $\ER^d$-valued fBm 
on $(-\infty,0]$.

Then, such an initial distribution is classically called an invariant distribution if it is 
invariant by the transformation ${\cal Q}_t$ for every $t\ge0$. However, the concept 
of uniqueness of invariant distribution is slightly different from the classical setting. 
Actually, if ${\bar Q}\mu$ stands for the distribution of the whole process $(X_t^\mu)_{t\ge0}$ 
with initial distribution $\mu$, one says that uniqueness of the invariant distribution holds 
if the stationary regime is unique (in other words, this concept of uniqueness 
corresponds to the classical one conditioned by the equivalence relation: 
$\mu\sim\nu \Longleftrightarrow \bar{Q}\mu\sim\bar{Q}\nu$, see \cite{hairer2} for 
background). We refer to \cite{hairer,hairer2,hairer-pillai} for criteria 
of uniqueness in different settings: additive noise, multiplicative noise with $H>1/2$ 
and multiplicative noise with $H\in(1/3,1/2)$ in the last one (in an hypoelliptic context).

The additive result of \cite{hairer} is obtained by a coupling strategy that we briefly 
recall here. Classically, coupling two paths issued of $\mu_0$ and $\mu$ where 
the second one denotes an invariant distribution of $({\cal Q}_t)_{t\ge0}$ consists 
in finding a stopping time ${\tau_{\infty}}$ such that 
$(X_{t+{\tau_{\infty}}}^{\mu_0})_{t\ge0}= (X_{t+{\tau_{\infty}}}^{\mu})_{t\ge0}$ 
(so that the rate of convergence in total variation can be derived from some 
bounds on  $\PE({\tau_{\infty}}>t)$, $t\ge0$). Now, let us detail the strategy. 
First, one classically waits that the paths get close. Then, at each trial, the 
coupling attempt is divided in two steps. First, one tries in Step 1 to stick the 
positions on an interval of length $1$. Then, in Step 2, one tries to ensure that 
the paths stay stuck until $+\infty$.  Actually, oppositely to the Markovian case 
where the paths stay naturally together after a clustering (by putting the same 
noise on each coordinate), the main difficulty here is that, due to the memory, 
staying together is costly. In other words, this property can be ensured only 
with the help of a  non trivial coupling of the noises. One thus talks of \textit{asymptotic 
coupling}. If one of the two previous steps fails, we will begin a new attempt but 
only after a (long) waiting time which is called Step $3$. During this step, one 
again waits that the paths get close but one also expects the memory of the 
coupling cost to vanish sufficiently in order to begin the new trial with a weak 
weight of the memory.

In the previous construction, the fact that $\sigma$ is constant is fundamental 
for ensuring the two following properties:
\begin{itemize}
\item If two paths $B^1$ and $B^2$ of the fBm differ from a drift term, then 
two paths $X^1$ and $X^2$ of \eqref{fractionalSDE0} respectively directed by 
$B^1$ and $B^2$ also differ from a drift term, which allows in particular to 
use Girsanov Theorem to build the coupling in Step 1.   
\item Under some ``convexity'' assumptions on the drift apart from a 
compact set, two paths $X^1$ and $X^2$ directed by the same fBm 
(or more precisely, by two slightly different paths) get closer and the 
distance between the two paths can be controlled deterministically. 
\end{itemize}     

In the current work,  $\sigma$ is not constant and the two above properties 
are no longer valid. The challenge then is to  extend the applicability of the 
previous coupling scheme to such a situation. The replacement of  each of 
the above properties  requires us to deal with different (though related) 
difficulties.   In order to be able to extend the  Girsanov argument  used in 
Step 1 to a non constant $\sigma$, we will restrain ourselves to diffusion 
coefficients  for which some injective function  of two copies of the process 
differs by a drift term whenever their driving fBm do. A natural assumption 
on $\sigma$ granting the latter property is  that $x\mapsto\sigma^{-1}(x)$ 
is (well-defined and is) a Jacobian matrix. This will be the setting of the 
present paper.

 As concerns a suitable substitution  of the second lacking property, a natural (but 
 to our knowledge so far not explored) idea is  to try to extend Meyn-Tweedie 
 techniques (see $e.g.$ \cite{DownMeynTweedie} for background) to the fractional 
 setting. More precisely, even if the paths do not get closer to each other 
 deterministically, one could expect that  some Lyapunov assumption could 
 eventually  make the two paths return in some compact set simultaneously. 
 The main contribution of the present paper is to incorporate such a Lyapunov-type 
 approach into the study of long-time convergence  in the fractional diffusion setting. 
 As one could expect, compared to the Markovian case, the problem is much more 
 involved. Actually, the  return  time to a compact set after a (failed) coupling attempt 
 does not only depend on   the positions of the processes after it, but also on all the 
 past of the fBm. Therefore,  in order that  the coupling  attempt succeeds  with 
 lower-bounded probability, one needs to establish some controls on the past 
 behavior of the fBms  that drive the two copies of the process,  conditionally to 
 the failure of the previous attempts. This point is one of the main difficulties of 
 the paper,  since, in the corresponding estimates, we  carefully have to take 
 into account   all the deformations  of the distribution that  previously failed 
 attempts induce. Then, we show that after a sufficiently long waiting time, 
 conditionally on previous fails the probability that the two paths be in a compact 
 set and  that the influence  of past noise on the future be controlled, is lower-bounded.  
 Bringing all the estimates together yields a global control of the coupling time and a rate of convergence which is similar to the one in \cite{hairer} in the additive noise case.
 
\subsection{Assumptions and Main Result}

Remind that in the whole section it is assumed that $H\in (1/2,1)$.  
We begin by a condition for the existence and uniqueness of solutions for \eqref{fractionalSDE0}:

\noindent $\mathbf{(H_0)}$: $b$ is a locally Lipschitz and sublinear function and 
$\sigma$ is a bounded $(1+\gamma)$-Lipschitz continuous function with 
$\gamma\in(\frac{1}{H}-1,1]$ ($i.e.$ $\sigma$ is a ${\cal C}^1$-function whose 
partial derivatives are bounded and globally $\gamma$-Hölder-continuous).\smallskip

\noindent Before introducing the second assumption, let us give a definition. We say 
that a function $V:\ER^d\rightarrow\ER$ is {\em essentially quadratic} if  it is a positive 
${\cal C}^1$-function such that $\nabla V$ is Lipschitz continuous and such that
\[
 \liminf_{\vert x\vert\to\infty} \frac{V(x)}{|x|^2}>0 \quad \mbox{  and }\quad 
  \abs[\nabla V]\le C \sqrt{V}\quad (C>0).
 \]
In order to ensure the existence of the invariant distribution, we now introduce a 
Lyapunov-stability assumption $\mathbf{(H_1)}$ through such a function $V$: 
 
\noindent $\mathbf{(H_1)}$: There exists an essentially quadratic function 
$V:\ER^d\rightarrow\ER$, there exist some positive $\beta_0$ and $\kappa_0$ such that
\[
\forall\;x\;\in\ER^d,\quad (\nabla V(x)|b(x))\le \beta_0-\kappa_0 V(x).
\]

\begin{rmrk}[Comparison to the Markovian case]
For the coupling strategy, the above assumption will be certainly used to 
ensure that the paths live in a compact set with a high probability. Note that in 
the classical diffusion setting, such a property holds with some less restrictive 
Lyapunov assumptions. Here, the assumptions $\mathbf{(H_0)}$ and $\mathbf{(H_1)}$  
essentially allow us to consider only (attractive) drift terms whose growth is linear 
at infinity.  This  more restrictive condition can be understood as a consequence 
of the lack of martingale property for the integrals driven by fBms, which leads 
in fact to a more important contribution of the noise component.
\end{rmrk} 

\begin{rmrk}[Comparison to previous results]
In \cite{hairer}, the  corresponding assumption is a contraction condition 
out of a compact set: for any $x,y$, $(b(x)-b(y)|x-y)\le \beta_0-\kappa_0|x-y|^2$. 
This means that even in the constant case, our work can cover some new cases. 
For instance, if $d=2$ and $b(z)=-z-\rho \cos(\theta_z) z^\perp$ 
(where $\rho\in\ER$, $\theta_z$ is the angle of $z$ and $z^\perp$ is its normal vector), 
Assumption $\mathbf{(H_1)}$ holds whereas one can check that the contraction 
condition is not satisfied if $\rho>2$. 
\end{rmrk} 

When the paths are in this compact set, one tries classically to couple them with 
positive probability. But, as mentioned before, the specificity of the non-Markovian 
setting is that the coupling attempt generates a cost for the future (in a sense made 
precise later). In order to control this cost or more precisely in order to couple the 
paths with the help of a controlled drift term, we need to ensure the next assumption: 
 
\smallskip
\noindent $\mathbf{(H_2)}$ $\forall x\in\ER^d$, $\sigma(x)$ is invertible and there 
exists a ${\cal C}^1$-function $h=(h_1,\ldots,h_d):\ER^d\rightarrow\ER^d$ such 
that the Jacobian matrix $\nabla h=(\partial_{x_j} h_i)_{i,j\in\{1,\ldots,d\}}$ satisfies 
$\nabla h(x)=\sigma^{-1}(x)$ and such that $\nabla h$ is a locally Lipschitz function on $\ER^d$.

\smallskip
\noindent
\begin{rmrk}[On the regularity of the diffusion matrix]
Under $\mathbf{(H_0)}$ and $\mathbf{(H_2)}$, $h$ is a global 
${\cal C}^1$-diffeomorphism from $\ER^d$ to $\ER^d$. Indeed, under these assumptions,  
$\nabla h$ is invertible everywhere and $x\mapsto [(\nabla h)(x)]^{-1}=\sigma(x)$ 
is bounded on $\ER^d$. Then, the property (which will be important in the sequel), 
follows from the Hadamard-L\'evy theorem (see e.g. \cite{ruadulescu}). 

As mentioned before, the main restriction here is to assume that 
$x\mapsto \sigma^{-1}(x)$ is a Jacobian matrix. However, note that there is no 
assumption on $h$ (excepted smoothness). In particular, $\sigma^{-1}$ does not 
need to bounded. This allows us to consider for instance some cases where 
$\sigma$ vanishes at infinity.
\end{rmrk}

\smallskip
\noindent Let us exhibit some simple classes of SDEs  for which $\mathbf{(H_2)}$ 
is fulfilled. First, it contains the class of non-degenerated SDEs for which each 
coordinate is directed by one real-valued fBm. More precisely,  if for every 
$i\in\{1,\ldots,d\}$, 
\[
dX_t^i=b_i(X_t^1,\ldots,X_t^d)dt+\sigma_i(X_t^1,\ldots,X_t^d) dB_t^i
\]
where $\sigma_i:\ER^d\rightarrow\ER$ is a   ${\cal C}^1$ positive function, 
Assumption~$\mathbf{(H_2)}$ holds. Now, let us also remark that since, for a 
given constant matrix, $\nabla (P h)=P\nabla h$, we have the following equivalence:
\[
\textnormal{$\exists$ $h$ such that $\nabla h=\sigma^{-1}$ 
$\Longleftrightarrow$ 
$\exists$ $\tilde{h}$, $\exists$ an invertible matrix $P$ such that $\sigma^{-1}=P\nabla \tilde{h}$}.
\]
One deduces from this property that  $\mathbf{(H_2)}$ also holds true if:
\bqne
\sigma (x)=P {\rm Diag}\left(\sigma_1(x_1,\ldots,x_d),\ldots, \sigma_d(x_1,\ldots,x_d)\right)
\eqne
where $P$ is a given invertible $d\times d$-matrix and for every $i\in\{1,\ldots,d\}$ 
$\sigma_i$ has the same properties as before.

We are now able to state our main result. One denotes by 
${\cal L}((X_{t}^{\mu_0})_{t\ge0})$ the distribution of the process on 
the set ${\cal C}([0,+\infty),\ER^d)$ starting from an initial distribution $\mu_0$
 and by $\bar{\cal Q}\mu$ the distribution of the stationary solution (starting 
 from an invariant distribution $\mu$). The distribution $\bar{\mu}_0(dx)$ 
 denotes the first marginal of $\mu_0(dx,dw)$.

\begin{thrm}\label{theo:principal} Let $H\in(1/2,1)$. Assume 
$\mathbf{(H_0)}$, $\mathbf{(H_1)}$ and $\mathbf{(H_2)}$. Then, existence 
and uniqueness hold for the invariant distribution $\mu$ (up to equivalence). 
Furthermore, for every initial distribution $\mu_0$ for which there exists $r>0$ 
such that { $\int |x|^r \bar{\mu}_0(dx)<\infty$}, for each $\varepsilon>0$ there 
exists $C_{\varepsilon}>0$ such that 
\[
\| {\cal L}((X_{t+s}^{\mu_0})_{s\ge0})-\bar{\cal Q}\mu\|_\mathrm{TV}
\leq C_\varepsilon t^{-(\frac{1}{8}-\varepsilon)}.
\]
\end{thrm}

\begin{rmrk}
In the previous result, the main contribution is the fact that one is able to 
recover the rates of the additive case. Existence and uniqueness results are not 
really new. However, compared with the assumptions of \cite{hairer2}, one observes 
that when  $x\mapsto\sigma^{-1}(x)$ is a Jacobian matrix (assumption which does 
not appear in \cite{hairer2}), our other assumptions are slightly less constraining. 
In particular, $b$ is assumed to be locally Lipschitz and sublinear (instead of 
Lipschitz continuous) and, as mentioned before, $x\mapsto \sigma^{-1}(x)$ does 
not need to bounded.
\end{rmrk}

\noindent \textbf{Some ingredients of the proof.} The aim of this section is 
to give some ideas of the proof. As explained above, the scheme is similar 
to that of \cite{hairer}. The starting point is to consider a couple $(X,\tilde{X})$ 
of solutions to \eqref{fractionalSDE0} with respective driving fBms denoted by 
$B$ and $\tilde{B}$. The underlying innovation processes are denoted by $W$ 
and $\tilde{W}$. For the sake of simplicity, assume that $X_0=x\in\ER^d$ and 
that the initial condition of  $\tilde{X}$ is the invariant distribution $\mu$. For 
every $k\ge1$, denote by $\tau_{k-1}$, the starting time of the  $k^{th}$-coupling 
attempt and by $\Delta \tau_k$ its  duration. If the coupling is successful, $X$ 
and $\tilde{X}$ get stuck after Step 1, $i.e.$ $\tau_{k-1}+1$. We thus define 
$\tau_\infty:=\tau_{k^*-1}+1$ where $k^*:=\inf\{k\ge1,\Delta \tau_k=+\infty\}$. By 
construction, 
\[
\forall t\ge0,\quad  \| {\cal L}((X_{t+s}^{\mu_0})_{s\ge0})-\bar{\cal Q}\mu\|_\mathrm{TV}\le \PE(\tau_\infty>t).
\]
Furthermore, using that 
\[
\forall t\ge 1,\quad \PE(\tau_\infty>t-1)=\PE(\tau_0+\sum_{k=1}^{+\infty}\Delta \tau_k 1_{k^*>k}>t),
\]
it can be shown (see \cite{fontbona-panloup}, section 5) that the problem can be 
more or less reduced to the control (uniform  in $k$) of : 
\[
\quad \PE(\Delta \tau_k<+\infty|\Delta \tau_{k-1}<+\infty,{\cal F}_{\tau_{k-1}})
\quad \textnormal{and}\quad 
\ES[|\Delta \tau_k|^p|\Delta \tau_k,{\cal F}_{\tau_{k-1}}]
\]
where $({\cal F}_t)_{t\ge0}$ is the usual augmentation of  
$((\sigma(W_s,\tilde{W}_s))_{s\le t})_{t\ge0}$ and $p\in(0,1)$ is a real that 
one will try to maximize (as suggested by Theorem \ref{theo:principal}, this 
expectation will be finite as soon as $p<1/8$).

\noindent \textbf{$(K,\alpha)$-admissibility.} The quantities defined previously 
can be controlled only if the positions of each component and the past of their 
noise satisfy some conditions at the beginning of the attempt. One will talk 
about $(K,\alpha)$-admissibility. In order to define this concept, we now assume that 
$W$ and $\tilde{W}$ differ from a drift term denoted by $g_w$ : 
$dW_t=d\tilde{W}_t+g_w(t)dt$. The function $g_w$ will be supposed to be 
null before $\tau_0$, $i.e.$ before the first attempt and also during Step $3$. 
In order to quantify the impact of $g_w$ on the future attempts, one introduces 
an operator  ${\cal R}_T$ defined (when it makes sense) by: for all $T\ge0$ 
and for all  $g:\ER\rightarrow\ER$ 
\[
({\cal R}_T g)(t)=\int_{-\infty}^0 \frac{t^{\frac{1}{2}-H} (T-s)^{H-\frac{1}{2}}}{t+T-s} g(s) ds,
\quad t\in(0,+\infty).
\]
The admissibility condition can be then defined as follows:
\begin{dfntn} Let $K$ and $\alpha$ be some positive numbers and $\tau$ 
be a  $({\cal F}_t)_{t\ge0}$-stopping time.  One says that the system is 
$(K,\alpha)$-admissible at time $\tau$  if $\tau(\omega)<+\infty$ and if 
$(X_\tau^{1}(\omega),X_\tau^{2}(\omega), (W^1(\omega),W^2(\omega))_{t\le \tau})$ 
satisfies :
\begin{equation}\label{hairerassumpcond}
 \sup_{T\ge0} \int_0^{+\infty} (1+t)^{2\alpha}|({\cal R}_T g_w^\tau)(t)|^2dt
 \le 1 \quad \text{where } g_w^\tau(.)=g_w(.+\tau)
\end{equation}
is the shift of drift term between $W$ and $\tilde W$, 
and if $(X_\tau^{1}(\omega),X_\tau^{2}(\omega), (W^1(\omega),W^2(\omega))_{t\le \tau})
\in \Omega_{K,\alpha,\tau}$ where
\begin{equation}\label{Kadmiscond}
\Omega_{K,\alpha,\tau}:=
\{ |X_\tau^{1}(\omega)|\le K,\; |X_\tau^{2}(\omega)|\le K, 
\;  \varphi_{\tau,\varepsilon_\theta}(W^1(\omega))\le K \;
\quad\textnormal{and}\quad
\; \varphi_{\tau,\varepsilon_\theta}(W^2(\omega))\le K\},
\end{equation}
and  $\varepsilon_\theta=\frac{H-\theta}{2}$ with $\theta\in(1/2,H)$ and for all  $\varepsilon>0$,
\[
\varphi_{\tau,\varepsilon}(w)=
\sup_{\tau\le s\le t\le \tau+1}\left\vert 
\frac{1}{t-s}\int_{-\infty}^{\tau-1} (t-r)^{H-\frac{1}{2}}-(s-r)^{H-\frac{1}{2}} dw_r\right\vert 
+\sup_{\tau-1\le u<v\le \tau}\frac{ |w(v)-w(u)|}{|v-u|^{\frac{1}{2}-\varepsilon}}.
\]
\end{dfntn}
The $(K,\alpha)$-admissibility is a sufficient condition to ensure that the coupling 
succeeds with lower-bounded probability ($i.e.$ such that the paths remain stuck 
until infinity). More precisely, under this condition, one is able to show that there 
exists $\delta_0>0$ such that for every $k\ge1$,
\[
\PE(\Delta \tau_k=+\infty|\{\Delta_{\tau_{k-1}}<+\infty\}\cap \Omega_{K,\alpha,\tau_{k-1}})
\ge\delta_0.
\]
Condition \eqref{hairerassumpcond}, which comes from \cite{hairer}, is adapted for 
Step $2$, ensuring that if $X_{\tau_{k-1}+1}=\tilde{X}_{\tau_{k-1}+1}$ ($i.e.$ if Step 
$1$ succeeds), one can build some couples $(W,\tilde{W})$ on successive intervals 
of lengths $2^{\alpha N}$ such that the probability that the paths remain stuck until 
infinity is lower-bounded. The second condition plays an important role for Step 1. 
When $\omega\in \Omega_{K,\alpha,\tau_{k-1}}$, the cost to stick the paths can  
be uniformly controlled under Assumption $\mathbf{(H_2)}$. Then, one of the main 
difficulties of the proof in the multiplicative case is to  show that the probability of 
$(K,\alpha)$-admissibility is also lower-bounded: one needs to show that there exists 
$\delta_1>0$ such that for every $k\ge1$,
\[
\PE(\Omega_{K,\alpha,\tau_{k}}|\Delta \tau_{k-1}<+\infty)\ge\delta_1.
\]
The proof of this property is achieved in two main steps. In the first one, one needs 
to control the  increments of $W$ and $\tilde{W}$ before $\tau_{k}$, conditionally 
to the failure of the previous attempts. This control requires a sharp knowledge of 
the construction of the innovations used in Steps 1 and 2 in order to identify the 
 distortions  generated by each scenario of failed attempt. Then, plugging these 
 controls into the Mandelbrot-Van Ness representation allows us to estimate the 
 impact of these distortions on the future of the increments of  the fBm (see 
 Lemmas 4.5, 4.6 et 4.7 of \cite{fontbona-panloup} for more details).
 
With the help of the Lyapunov assumption $\mathbf{(H_1)}$, one tries in the 
second step to upper-bound the quantity  
$\ES[V^r(X_{\tau_k})+V^r(\tilde{X}_{\tau_k})|\Delta \tau_k<+\infty]$ (for a 
positive power $r$). Conditionally to the first step, the keypoint is the  
following contraction property (Proposition 4.4 of \cite{fontbona-panloup}) : 
if $\mathbf{(H_0)}$ and $\mathbf{(H_1)}$ hold, then,  for every $\theta\in (1/2,H)$, 
there exists $\bar\rho\in(0,1)$, $C>0$ and $r>0$ such that for every starting point $x\in\ER^d$,
\begin{equation}\label{eq:lyap1}
 V^r(X_1)\le \bar\rho V^r(x)+C(1+\|B^H\|_\theta^{0,1})
 \quad\textnormal{where}\quad 
 \|B\|_\theta^{0,1}=\sup_{0<s<t<1}\frac{|B^H_t-B^H_s|}{(t-s)^{\theta}}.
\end{equation}
We refer to \cite{fontbona-panloup} for the complete proof.

\section{The penalized bandit process}\label{Sec_PenalizedBandit}


The \emph{two-armed bandit algorithm} is a theoretical procedure to choose 
asymptotically the most profitable arm of a slot machine, or bandit; it was also 
used in the fields of mathematical psychology and of engineering. This algorithm 
has been widely studied, for instance in \cite{LPT04,LP08b}. The key idea is to 
use a (deterministic) sequence of learning rates, rewarding an arm if it delivers 
a gain. Depending on the speed of convergence to 0 of this sequence, the 
algorithm is often faillible (it would not always select asymptotically the right arm, 
see \cite{LPT04}).

It is possible to improve its results and ensure infaillibility by introducing penalties 
when the arm does not deliver a gain: this modification is called the 
\emph{penalized bandit algorithm}, and it is studied in \cite{LP08}. The authors 
show that, with a correct choice of penalties and rewards, and with the appropriate 
renormalization, the algorithm converges weakly to a probability measure $\pi$, 
which is the stationary distribution of the Piecewise Deterministic Markov Process 
with following infinitesimal generator
\[
\mathcal Lf(x)=(1-p-px)f'(x)+qx\frac{f(x+g)-f(x)}{g}.
\]
where $0<q<p<1$ being the respective probabilities of gain of the 
two arms and the positive parameter $g$ runs the asymptotic behaviour of the sequences of 
the rewards and penalties (see Section~3 in \cite{LP08} for more details). This process is 
also studied in~\cite{GPS}. For the 
sake of simplicity, we set $g=1$ in the sequel. Moreover, since the interval $[0,(1-p)/p)$ 
is transient, computations are easier if we study the translated process 
$Y=X-\frac{1-p}p$, driven by the following generator:
\begin{equation}\label{Eq_GenBandit}
\mathcal L^Yf(y)=-pyf'(y)+q\left(y+\frac{1-p}p\right)\big(f(y+1)-f(y)\big).
\end{equation}
It is possible to deduce the dynamics of the process from the generator 
(see \cite{Dav93}): between the jumps, $Y$ evolves as the solution of the ODE $y'_t=-py_t$, 
and it jumps with jump rate $t\mapsto\zeta(Y_t)=q\left(Y_t+\frac{1-p}p\right)$ from 
$Y_t$ to $Y_t+1$.

In \cite{LP08}, the authors show that $\pi$ admits a density with support $[(1-p)/p,+\infty)$ 
and exponential moments of order up to $u_M$, where $u_M$ is the unique positive solution 
of the equation
\begin{equation}
\frac{\exp(u_M)-1}{u_M}=\frac pq.
\label{Eq_uM}\end{equation}
We shall prove the latter too, but with a different argument (see Remark~\ref{Rk_LaplTrans}).

In the sequel, we denote by \emph{penalized bandit process} the process with 
initial distribution $\mu_0$ following the dynamics of $\mathcal L^Y$, and by $\mu_t$ 
its law at time $t\geq 0$. This section is devoted to the long time behavior 
of this process with respect to Wasserstein and total variation distances. The 
estimates rely on the construction of explicit couplings. This approach is closely 
related to the paper~ \cite{Bou14} which is dedicated to the study of a Piecewise 
Deterministic Markov Process related to a pharmacokinetic model introduced in~\cite{BCT08}. 

\subsection{Wasserstein convergence}

Firstly, let us recall the definitions of Wasserstein and total variation distances 
between two measures:
\begin{align*}
\Wass_n(\mu,\nu)
&=\inf\left\{\E[|X-Y|^n]^{\frac1n}:(X,Y)\text{ coupling of }\mu\text{ and }\nu\right\},\\
\|\mu-\nu\|_{\mathrm{TV}}
&=\inf\left\{\prob(X\neq Y):(X,Y)\text{ coupling of }\mu\text{ and }\nu\right\}.
\end{align*}
In the following, let $\mu_0$ and $\widetilde\mu_0$ be two probabilities on $\R_+$. 
The following proposition holds:
\begin{prpstn}
\label{WassBandit}
We have, for all $t\geq0$
\begin{equation}
\label{Eq_WassBandit}
\Wass_1(\mu_t,\widetilde\mu_t)\leq\Wass_1(\mu_0,\widetilde\mu_0)\e^{-(p-q)t}.
\end{equation}
\end{prpstn}

\begin{proof}
Let $(Y,\widetilde Y)$ be generated by
\begin{align}\label{Eq_GenBandit2}
\mathcal L_2^Yf(y,\widetilde y)=
&-py\partial_yf(y,\widetilde y)-p\widetilde y\partial_{\widetilde y}f(y,\widetilde y)\notag\\
		&+q(y-\widetilde y)(f(y+1,\widetilde y)-f(y,\widetilde y))
		+q\left(\widetilde y+\frac{1-p}p\right)(f(y+1,\widetilde y+1)-f(y,\widetilde y)),
\end{align}
for $y\geq\widetilde y$, and of symetric expression for $\widetilde y\geq y$, 
and such that $(Y_0,\widetilde Y_0)$ is a coupling of $(\mu_0,\widetilde\mu_0)$ 
realizing $\Wass_1(\mu_0,\widetilde\mu_0)=\E\left[|Y_0-\widetilde Y_0|\right]$. 
It is easy to check that \eqref{Eq_GenBandit2} reduces to \eqref{Eq_GenBandit} 
if $f$ depends only on $y$ or $\widetilde y$, which means that 
$(Y_t,\widetilde Y_t)_{t\geq0}$ generated with $\mathcal L_2^Y$ is a coupling 
of $(\mu_t,\widetilde\mu_t)_{t\geq0}$. With this coupling, either the higher process 
jumps alone or the jump is simultaneous. It is easy to check that this coupling is 
monotonous, \emph{i.e.} for all $t\geq0,(Y_t-\widetilde Y_t)(Y_0-\widetilde Y_0)\geq0$. 
Monotonicity comes from the fact that the higher process jumps more often but 
stays above the other since the jumps are positive. Assume that 
$Y_0\geq\widetilde Y_0$. By monotonicity, we have, for all $t\geq0$,
\[
\E\left[|Y_t-\widetilde Y_t|\right]=\E[Y_t]-\E[\widetilde Y_t],
\]
so all we have to do is to study $h:t\mapsto\E[Y_t]$. With $f(y)=y$, \eqref{Eq_GenBandit} 
leads to $\mathcal Lf(y)=\frac{q(1-p)}p-(p-q)y$, and, by Dynkin formula, the function $h$ 
satisfies the ordinary differential equation $h'(t)=\frac{q(1-p)}p-(p-q)h(t)$.
One deduces immediately that
\begin{equation}\label{Eq_Wass1}
\E[Y_t]=\frac{q(1-p)}{p(p-q)}+\left(\E[Y_0]-\frac{q(1-p)}{p(p-q)}\right)\e^{-(p-q)t},
\end{equation}
(recall that $p>q$). Then,
\[\E\left[|Y_t-\widetilde Y_t|\right]=\E\left[Y_0-\widetilde Y_0\right]\e^{-(p-q)t},\]
which leads directly to \eqref{Eq_WassBandit}
\end{proof}

\begin{rmrk}\label{Rk_LaplTrans}
The Dynkin formula is a powerful tool for studying the moments of Markov processes. 
One can use it with $f(y)=\e^{uy}$ to study the Laplace transform of the transient 
process $\psi(t,u)=\E[\e^{uY_t}]$. We have
\[
\mathcal L^Yf(y)=q\frac{1-p}p(e^u-1)f(y)+(q(e^u-1)-up)yf(y),
\]
so $\psi$ satisfies the following PDE:
\[\partial_t\psi(t,u)=q\frac{1-p}p(e^u-1)\psi(t,u)+(q(e^u-1)-up)\partial_u\psi(t,u).\]
If $\mu_0$ is the invariant measure $\pi$, then $\partial_t\psi(t,u)=0$, and the Laplace 
transform $u\mapsto \psi(u)$ is solution of the following ODE 
\[
\partial_u\left(\log(\psi(u))\right)=\frac{q\frac{1-p}p(\e^u-1)}{up-q(\e^u-1)},
\]
and the right-hand side is finite for $u\in[0,u_M)$, when $u_M$ is the solution of 
Equation~\eqref{Eq_uM}.
\end{rmrk}

Note that the set of polynomials of degree $n$ is stable under the action of $\mathcal L^Y$. 
This is an important property, since it theoretically enables us to compute the moments of 
$Y_t$ by induction, with the help of Dynkin formula, just as we did for the first moment in 
the proof of Proposition~\ref{WassBandit}. Similarly, it is possible to study the function 
$h_n(t)=\E[|Y_t-\widetilde Y_t|^n]$ which provides an upper bound of 
$\Wass_n(\mu_t,\widetilde\mu_t)$. Indeed, we have, for $f(y,\widetilde y)=|y-\widetilde y|^n$,
\[
\mathcal L^Y_2f(y,\widetilde y)=
-n(p-q)|y-\widetilde y|^n+q\sum_{k=0}^{n-2}{\binom nk|y-\widetilde y|^{k+1}},
\]
so
\[
h'_n(t)=-n(p-q)h_n(t)+q\sum_{k=0}^{n-2}{\binom nkh_{k+1}(t)}.
\]
Then, using Gronwall lemma easily leads to check, by induction, that 
$h_n(t)=\mathcal{O}(\e^{-n(p-q)t})$. Which leads to the following result:

\begin{prpstn}\label{WassBandit2}
For all $n\in\N^\star$, there exists a positive constant $C$ such that, for all $t\geq0$,
\[
\Wass_n(\mu_t,\widetilde\mu_t)\leq C\e^{-(p-q)t}.
\]
\end{prpstn}

\subsection{Total variation convergence}

\begin{figure}[!ht]
\centering
\begin{tikzpicture}
\draw[blue] (0,2) node[left]{$\widetilde Y_0$};
\draw[domain=0:1,smooth,blue] plot (\x,{2*exp(-0.5*\x)});
\draw[dashed,blue] (1,1.21306) -- (1,2.21306) node[above]{$\widetilde Y_{\widetilde T}$};
\draw[domain=1:1.948154,smooth,blue] plot (\x,{2.21306*exp(-0.5*(\x-1))});
\draw[red] (0,1) node[left]{$Y_0$};
\draw[domain=0:1.948154,smooth,red] plot (\x,{exp(-0.5*\x)});
\draw[dashed,red] (1.948154,0.3775407) -- (1.948154,1.3775407) node[above]{$Y_T$};
\draw[domain=1.948154:2.5,smooth,purple] plot (\x,{1.3775407*exp(-0.5*(\x-1.948154))});
\end{tikzpicture}
\caption{Expected behaviour of the coalescent coupling for the penalized bandit process.}
\label{Fig_TVBandit}
\end{figure}
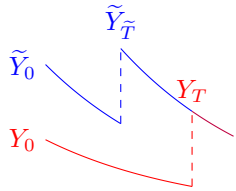

In the case of the penalized bandit process, total variation convergence is slightly 
harder than in \cite{Bou14}, since the jumps are always of amplitude range 1. Instead, 
we are going to use the arguments introduced in \cite{BCGMZ13}, based on the 
following observation: if $Y$ and $\widetilde Y$ are close enough, we can make 
them jump, not simultaneously like before, but with a slight delay for one of the copies, 
which would make it jump on the other one, as illustrated in Figure~\ref{Fig_TVBandit}.

In the following, denote by $\tau=\inf\{t\geq0:\forall s\geq0,Y_{t+s}=\widetilde Y_{t+s}\}$ 
the coalescence time of $Y$ and $\widetilde Y$. The goal of the sequel is to obtain 
exponential moments for $\tau$ (which may happen for correct couplings) and then 
use the total variation classic coupling inequality:
\[
\|\mu_t-\widetilde\mu_t\|_\mathrm{TV}\leq\prob(Y_t\neq\widetilde Y_t)\leq\prob(\tau>t).
\]
We have the following lemma:
\begin{lmm}\label{TVBandit1}\
Assume there exist positive constants $\bar y<+\infty,\bar\eps<1$ such that 
$Y_0,\widetilde Y_0\leq \bar y$ and $|Y_0-\widetilde Y_0|\leq\bar\eps$. Then, 
there exist a coupling $(Y_t,\widetilde Y_t)_{t\geq0}$ of $(\mu_t,\widetilde \mu_t)_{t\geq0}$ 
and an explicit positive constant $C(\bar y,\bar\eps)<+\infty$ such that, for all $t>0$, 
\begin{equation}
\label{Eq_TVBandit1}
\prob(\tau>t)\leq C(\bar y,\bar\eps)\left(\exp\left(-\frac{q(1-p)}pt\right)+\bar\eps\right).
\end{equation}
\end{lmm}

\begin{proof}
First, assume that $Y_0$ and $\widetilde Y_0$ are deterministic, and denote by 
$y=Y_0,\eps=\widetilde Y_0-y$. We assume w.l.o.g. that $\eps>0$. Let $T$ (resp. 
$\widetilde T$) be the first jump time of the process $Y$ (resp. $\widetilde Y$). 
Following the heuristics of Figure~\ref{Fig_TVBandit}, it is straightforward that
\[
Y_T=\widetilde Y_T\Leftrightarrow T=\frac1p\log\left(\eps+\e^{p\widetilde T}\right).
\]
Easy computations lead to
\[
\prob\left(\frac1p\log\left(\eps+\e^{p\widetilde T}\right)\leq s\right)
=\prob\left(\widetilde T\leq\frac1p\log\left(\e^{ps}-\eps\right)\right)=1-\Phi_y(s,\eps),
\]
with
\[
\Phi_y(s,\eps)=\exp\left(-\frac qp\left(\frac{1-p}p\log(\e^{ps}-\eps)+(y+\eps)
\left(1-\frac1{\e^{ps}-\eps}\right)\right)\right).
\]
As a consequence, the random variables $T$ and $\frac1p\log(\eps+\exp(p\widetilde T))$ 
admit densities w.r.t. the Lebesgue measure, which are respectively $f_y(\cdot,0)$ 
and $f_y(\cdot,\eps)$, with, for all $s\geq0$,
\[
f_y(s,\eps)=
\frac{q\e^{ps}}p\left(\frac{1-p}{\e^{ps}-\eps}+\frac{p(y+\eps)}{(\e^{ps}-\eps)^2}\right)\Phi_y(s,\eps).
\]
Let $T$ and $\frac1p\log(\eps+\exp(p\widetilde T))$ follow the $\gamma$-coupling 
(the coupling minimizing the total variation of their laws, see \cite{Lin02}). It is not 
hard to deduce from the very construction of this coupling the following equality:
\[
\prob\left(T=\frac1p\log\left(\eps+\e^{p\widetilde T}\right),T<t\right)
=\int_0^t{f_y(s,0)\wedge f_y(s,\eps)ds},
\]
where $x\wedge y=\min(x,y)$ and then
\begin{align}
\prob\left(\tau\leq t\right)&=\prob\left(Y_t=\widetilde Y_t\right)\geq
\prob\left(T=\frac1p\log\left(\eps+\e^{p\widetilde T}\right),T<t\right)\notag\\
	&\geq1-\frac12\left(\Phi_y(t,0)+\Phi_y(t,\eps)+\int_0^t{|f_y(s,0)-f_y(s,\eps)|ds}\right).
	\label{Eq_TV1}
\end{align}
The following upper bound is easily obtained, for any $0\leq\eps\leq\bar\eps$ 
and any $0\leq y\leq \bar y$:
\begin{equation}\label{Eq_TV2}
\Phi_y(s,\eps)\leq C_1\exp\left(-\frac{q(1-p)}ps\right),
\end{equation}
with $C_1=\exp\left(\frac{-q}p\left(\frac{1-p}p\log(1-\bar\eps)
+(\bar y+\bar\eps)\left(1-\frac1{1-\bar\eps}\right)\right)\right)$. In order to apply 
the mean-value theorem, we differentiate $f_y$ with respect to $\eps$. 
After some computations, one can obtain the following upper bound:
\begin{equation}\label{Eq_TV3}
\left|\frac{\partial f_y}{\partial\eps}(s,\eps)\right|\leq C_1C_2\exp\left(-\frac{q(1-p)}ps\right),
\end{equation}
where
\begin{align*}
C_2=&\frac{q^2((1-p)(1-\bar\eps)+p(\bar y+\bar\eps))}{p^2(1-\bar\eps)^2}
\left(1\vee\left(\frac{(2-p)(1-\bar\eps)+\bar y+\bar\eps}{(1-\bar\eps)^2}\right)\right)\\
	&+\frac{q(1-\bar\eps+2p(\bar y+\bar\eps))}{p(1-\bar\eps)^3}.
\end{align*}
Then, we easily have
\begin{equation}\label{Eq_TV4}
\int_0^t{|f_y(s,0)-f_y(s,\eps)|ds}
\leq C_1C_2\eps\int_0^{+\infty}{\exp\left(-\frac{q(1-p)}ps\right)ds}
\leq\frac{pC_1C_2}{q(1-p)}\bar\eps.
\end{equation}
Combining Equations~\eqref{Eq_TV1}, \eqref{Eq_TV2}, \eqref{Eq_TV3} 
and \eqref{Eq_TV4}, and denoting by $C(\bar y,\bar\eps)=C_1+\frac{pC_1C_2}{2q(1-p)}$, 
we have
\[
\prob\left(\tau>t\right)\leq C(\bar y,\bar\eps)\left(\exp\left(-\frac{q(1-p)}pt\right)+\bar\eps\right),
\]
then \eqref{Eq_TVBandit1} is straightforward. This upper bound does not depend on $Y_0$ 
and $\widetilde Y_0$, so this result still holds for random starting points, provided that 
they belong to $[0,\bar y]$.
\end{proof}

Proposition~\ref{WassBandit} and Lemma~\ref{TVBandit1} are the main tools to prove 
the following result:

\begin{prpstn}
\label{TVBandit}
Let $t_0>0$. There exists an explicit positive constant $K<+\infty$ such that, for all $t\geq t_0$,
\begin{equation}
\label{Eq_TVBandit}
\| \mu_t-\widetilde\mu_t\|_\mathrm{TV}
\leq K\e^{-vt},\qquad\text{with }v=\frac{p-q}{2+\frac{p(p-q)}{q(1-p)}}.
\end{equation}
\end{prpstn}

\begin{proof}
Let $\alpha\in(0,1)$ and $u>0$. We use the coupling from Proposition~\ref{WassBandit} 
in the domain $[0,\alpha t]$ and the coupling from Lemma~\ref{TVBandit1} in the domain 
$[\alpha t,t]$. We set $\bar\eps=\e^{-ut}$ and $\bar y=\frac{q(1-p)}{p(p-q)}+1$, and have 
the following inequality:
\begin{equation}\label{Eq_TV5}
\prob(\tau\leq t)\geq\prob\left(|Y_{\alpha t}-\widetilde Y_{\alpha t}|
\leq\bar\eps,Y_{\alpha t}\vee\widetilde Y_{\alpha t}\leq \bar y\right) 
\prob\left(\tau\leq t\left||Y_{\alpha t}-\widetilde Y_{\alpha t}|
\leq\bar\eps,Y_{\alpha t}\vee\widetilde Y_{\alpha t}\leq \bar y\right.\right).
\end{equation}
On the one hand,
\begin{align*}
\prob\left(|Y_{\alpha t}-\widetilde Y_{\alpha t}|>\bar\eps
\text{ or }Y_{\alpha t}\vee\widetilde Y_{\alpha t}>\bar y\right)
&\leq\prob\left(|Y_{\alpha t}-\widetilde Y_{\alpha t}|>\bar\eps\right)
+\prob\left(Y_{\alpha t}>\bar y\right)+\prob\left(\widetilde Y_{\alpha t}>\bar y\right)\\
	&\leq\left(\frac{\Wass_1(\mu_0,\widetilde\mu_0)}{\bar\eps}
	+\frac{\E[Y_0]+\E[\widetilde Y_0]}{\bar y-\frac{q(1-p)}{p(p-q)}}\right)\exp(-\alpha(p-q)t)\\
	&\leq C_3\exp((u-\alpha(p-q))t),
\end{align*}
with $C_3=\left(\Wass_1(\mu_0,\widetilde\mu_0)+\E[Y_0+\widetilde Y_0]\right)$. 
On the other hand, let $C_4=\sup_{t\geq t_0}C(\bar y,\e^{-ut})$. The constant 
$C_4$ is finite and, from Lemma~\ref{TVBandit1},
\begin{align*}
\prob\left(\tau>t\left||Y_{\alpha t}-\widetilde Y_{\alpha t}|\leq
\bar\eps,Y_{\alpha t}\vee\widetilde Y_{\alpha t}\leq \bar y\right.\right)
& \leq C_4\left(\e^{-ut}+\exp\left(-\frac{q(1-p)(1-\alpha)}pt\right)\right).
\end{align*}
Now, \eqref{Eq_TV5} reduces to
\[\prob(\tau>t)\leq1-\left(1-C_3\exp((u-\alpha(p-q))t)\right)
\left(1-C_4\left(\e^{-ut}+\exp\left(-\frac{q(1-p)(1-\alpha)}pt\right)\right)\right).\]
We optimize the rate of convergence by setting
\[
\alpha=\frac1{1+\frac{p(p-q)}{2q(1-p)}},\qquad u=\frac{(p-q)\alpha}2=v\text{ as defined above}.
\]
Then, \eqref{Eq_TVBandit} holds with $K=C_3+2C_4$.
\end{proof}

\section{Long time synchronization of large populations of interacting noisy rotators}

\subsection{The model}

Synchronization phenomena are a subject widely studied in physics and natural sciences.
Synchronization can occur in several different situations, for example
in the case of interacting cardiac cells, neurons, metronomes... (see \cite{cf:pikovsky}
for numerous examples of synchronization phenomena).
To construct a mathematical model in which a synchronization phenomenon
occurs, one may consider a population of dynamical systems
that interact with each other, and may perturb these systems
with noise (to modelize the internal noise of each dynamical system,
or the noise given by the interaction of each system with the surrounding environment).
We will focus here on a model given by a population of $N$ interacting
noisy rotators. Each rotator is defined by a phase
$\varphi_j(t) \in \mathbb{T}=\mathbb{R}/2\pi \mathbb{Z}$, and the evolution
of these phases is given by the following system
of stochastic differential equations:
\begin{equation}\label{eds:rotators}
 \dd \varphi_j(t) = -\frac KN \sum_{i=1}^N \sin(\varphi_j(t)-\varphi_i(t))\dd t +\sigma \dd B_j(t)\, ,
\end{equation}
where $K\geq 0$ and $\sigma>0$ are two constant parameters, and $(B_j)_{j=1,\ldots,N}$ is a
family of standard independent Brownian motions.
The interaction is of mean field type: the rotator $\varphi_j$ interacts with all the other
rotators, and the interaction term is constituted
of the sum of the contributions given by each one of these other rotators, with the same
weight $K/N$. 
Remark that the model is invariant by rotation: if $(\varphi_j(t))_{j=1,\ldots,N}$ is a solution
of \eqref{eds:rotators}, then it is also the case of
$(\varphi_j(t)+c)_{j=1,\ldots,N}$ for all real constant $c$.

\medskip

This model is known in the physics
literature as mean field plane rotors model, and is a particular
case of the celebrated noisy Kuramoto model (when the disorder follows the trivial 
distribution $\delta_0$).
Of course, since the dynamics of each isolated system
is very simple in this model (a Brownian motion on a circle), its aim
is not to describe a {\sl real} phenomenon (to do this one would need some
complex isolated system, in higher dimension and with several parameters),
but to provide a simple framework in which one can study analytically
a synchronization phenomenon.

\medskip

In this model we can speak of synchronization if the rotators concentrate around 
some phase, which we will
call center of synchronization. This may happen if the interaction is strong enough with respect
to the noise (since this later one incites the rotators to move independently), or in other
words if $K$ is sufficiently large compared to $\sigma$.
But with a simple time change one can replace $\sigma$ by $1$,
and $K$ by $K/\sigma^2$. So the {\sl real} parameter of the model is $K/\sigma^2$, and
we will set $\sigma=1$ in the remaining for simplicity.

\medskip

We will focus on the behavior of model in the limit of infinite size of population.
We will first consider the evolution on time intervals $[0,T]$ independent from 
$N$, when $N$ goes to infinity. We will then study the behavior of the model on 
longer time intervals, by making a rescaling in time depending on $N$. The content 
of the first part is based on the works \cite{cf:BGP,cf:GPP}, while
the second part describes the result proved in \cite{cf:BGPoq}.

\medskip

\subsection{Large populations and fixed time intervals}

We consider in this section the evolution on 
intervals of the type $[0,T]$, with $T$ independent from $N$.
In this case, since the coefficients in \eqref{eds:rotators} are smooth,
we can apply to our model the classical results
of the well-known theory of propagation
of chaos \cite{cf:Gartner,cf:Sznitman}. Let us
consider the empirical measure $\mu_{N,t}$ associated
to the model, i.e. the $\mathcal{M}_1(\mathbb{T})$-valued process 
$\mu_{N,t} = \frac 1N \sum_{i=1}^N \delta_{\varphi_j(t)}$,
where $\mathcal{M}_1(\mathbb{T})$ denotes the space of probability measures
on $\mathbb{T}$. If the initial condition $\mu_{N,0}$ converges weekly
to some $p_0\in \mathcal{M}_1$, then
the process $\mu_{N,t}$ converges weakly on $[0,T]$ to the deterministic trajectory
on $\mathcal{M}_1$ solution of the following Fokker-Planck type partially differential
equation:
\begin{equation}\label{edp:rotators}
 \partial_t p_t(\theta)= \frac 12 \partial^2_\theta p_t(\theta)- \partial_\theta
 \big[p_t(\theta)J*p_t(\theta)\big]\, ,
\end{equation}
where $*$ denotes the convolution operator and $J(\theta)=-K\sin(\theta)$.
In this limit $p_t$ is the distribution of the infinite population of rotators on $\mathbb{T}$
at time $t$.
The mass is preserved by this evolution, and
due to the presence of the Laplacian term this PDE admits a unique solution for any
initial condition $p_0\in \mathcal{M}_1$, which admits
a smooth and positive density for all $t>0$ (that we will also denote $p_t(\theta)$ 
for simplicity), element of  
$C^\infty((0,T)\times \mathbb{T}, \mathbb{R})$ (see for example \cite{cf:BGP} for
a proof of this regularity result).
Remark that the invariance by rotation of the finite size model is conserved in the limit,
since one can easily check that if $p_t(\theta)$ is solution of \eqref{edp:rotators},
then it is also the case of $p_t(\theta-c)$ for all $c\in \mathbb{R}$.

\medskip

A pleasant property of this PDE is that one can compute all its stationary solutions
in a semi-explicit way: $q$ is a (probability) stationary solution for \eqref{edp:rotators}
if and only if $q$ can be expressed as
\begin{equation}
 q(\theta)=q_{\psi,r}(\theta):=
 \frac{e^{2Kr \cos(\theta-\psi)}}{\int_{\mathbb{T}}e^{2Kr \cos(\theta'-\psi)}\dd\theta'}\, ,
\end{equation}
for some $\psi\in\mathbb{R}$ and some solution $r\geq 0$
of the fixed point problem
\begin{equation}\label{fixed point}
 r=\Psi(2Kr)\, ,
\end{equation}
where the fixed point function $\Psi$
is known explicitly, and satisfies some nice properties
that allow us to determine the number of solutions of the fixed point problem
according to the value of $K$: $\Psi$
is strictly concave and
bounded by $1$ on $(0,\infty)$, and satisfies $\Psi(0)=0$ and
$\Psi'(0)=1/2$ (see \cite{cf:BGP} for more details).

First remark that the equality $\Psi(0)=0$ implies that $r=0$ is always solution to the 
fixed point problem \eqref{fixed point}, which means that the uniform probability 
on the circle $q(\theta)=\frac{1}{2\pi}$ is always a stationary solution.
This solution corresponds to a total absence of synchronization in the model, since
in that case the population of rotators is distributed uniformly on the
$\mathbb{T}$. Moreover, the fact that $\Psi'(0)=1/2$ and the strict concavity of $\Psi$
imply that if $K\leq 1$, $r=0$ is the only solution
of \eqref{fixed point}. When $K\leq 1$ the interaction is too weak
to allow the apparition of synchronized states.

\begin{figure}\label{fig:stat states}
 \centering
 \includegraphics[height=6.5cm]{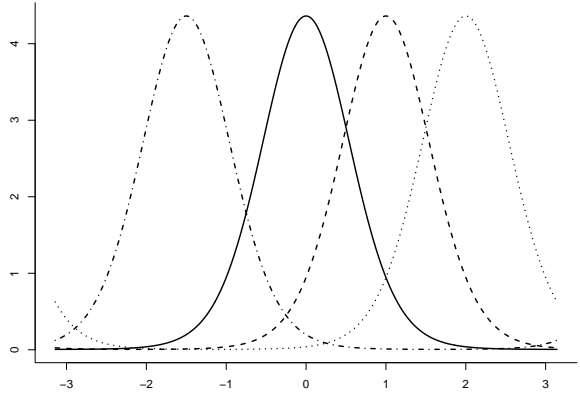}
 \caption{Graph of $q_{-1.5}$, $q_0$, $q_{1}$, $q_{2}$ (from the left to the right) when $K=2$.
 The elements of $M$ are the translations of the synchronized profile $q_0$.}
\end{figure}

On the other hand, if $K>1$, then there exists a unique positive solution $r_K$
to the fixed point problem, which means that in that case the set of stationary solutions
of \eqref{edp:rotators} is composed of $\frac{1}{2\pi}$ and 
a whole family $M$ of non-trivial stationary solutions, defined as
\begin{equation}
 M=\{q_\psi\, , \, \psi\in \mathbb{T}\}\, , 
 \quad \text{where} \quad q_\psi(\theta)= q_{\psi,r_K}(\theta)\, .
\end{equation}
Each stationary state $q_\psi$ is the translation by an angle $\psi$ of the profile
$q_0$, which corresponds to a concentration of the rotators around
the phase $0$ (see figure \ref{fig:stat states}).
So $q_\psi$ corresponds to a synchronization of the rotators around the
center of synchronization $\psi$.
From a geometrical point of view $M$ is a closed curve (in fact a circle since
all its points are the translation of the same profile) of synchronized stationary 
solutions, parametrized by their centers of synchronization.

\medskip

The next step in the understanding of the model is to determine the local stability of these different
stationary states. To do this let us linearize the evolution around these stationary states.
If we rewrite \eqref{edp:rotators} as $\partial_t p_t=F(p_t)$,
with $F(p)=\frac12 p''-(pJ*p)'$,
and consider a stationary state $q$ and
a smooth function $u$ satisfying 
$\int_\mathbb{T} u=0$, we can expand $F(q+u)$ as follows:
\begin{equation}
 F(q+u)=F(q)+\frac 12  u''(\theta)-
\big[ q(\theta)J*u(\theta)+ u(\theta)J* q(\theta)+ u(\theta)J*u(\theta)\big]'\, .
\end{equation}
$F(q)=0$ since $q$ is stationary, and by just keeping the linear terms we get 
a linearized evolution $\partial_t u=L_q u_t$ around the stationary profile $q$, with 
\begin{equation}\label{def Lq}
 L_q u(\theta) := \frac 12 u''(\theta)-
\big[ q(\theta)J*u(\theta)+ u(\theta)J* q(\theta)\big]'\, .
\end{equation}
The spectral properties of the operator $L_q$ determines the behavior of the 
solutions of \eqref{edp:rotators} in the neighborhood of $q$. 
Its spectrum can be obtained easily when $q=\frac{1}{2\pi}$ since in that case
$L_{\frac{1}{2\pi}}$ can be decomposed on the Fourier basis:
\eqref{def Lq} boils down to
\begin{equation}
 L_{\frac{1}{2\pi}} u(\theta) = \frac 12 u''(\theta)-\frac{1}{2\pi}\big[J*u(\theta)]'\, ,
\end{equation}
and by writing $u(\theta)
=\sum_{k=1}^\infty a_k \cos(k\theta)+\sum_{k=1}^\infty b_k\sin(k\theta)$ one simply gets
(recall that $J(\theta)=-K\sin\theta$):
\begin{equation}
 L_{\frac{1}{2\pi}}u(\theta)
 =-\frac 12(1-K)a_1\cos\theta-\frac12 (1-K)b_1\sin\theta-\frac12\sum_{k=2}^\infty
 k^2a_k\cos(k\theta)-\frac12\sum_{k=2}^\infty k^2 b_k\sin(k\theta)\, .
\end{equation}
This shows that when $K<1$ the spectrum of $L_{\frac{1}{2\pi}}$ is composed of 
strictly negative eigenvalues, and thus $\frac{1}{2\pi}$ is stable for the linearized evolution.
When $K>1$ the directions given by $\cos\theta$ and $\sin\theta$
become instable.

\medskip

The study of the spectrum of $L_{q_{\psi}}$ for $q_\psi\in M$ when $K>1$ is 
more involved, and requires the use of weighted Sobolev spaces. In \cite{cf:BGP} 
the authors show that in a well chosen weighted space the $L_{q_{\psi}}$ can also 
be decomposed on an orthogonal basis, and that the spectrum is made
of a decreasing sequence of non positive reals $(\lambda_i)_{i\geq 0}$ 
satisfying $\lambda_0=0$, $\lambda_i<0$ for
$i\geq 1$ and $\lambda_i\rightarrow_{i \rightarrow \infty} -\infty$.
The eigenvalue $0$ is associated to the eigenvector $q'_\psi$, which generates
the tangent space of $M$ at $q_\psi$.
So the fact that the operator $L_{q_{\psi}}$ has no effect in this direction
is not surprising, since in this direction the dynamics given by \eqref{edp:rotators}
is trivial. On the normal space (that is the subspace generated by the other eigenfunctions),
the linearized dynamics is stable with a spectral gap given by the largest negative
eigenvalue $\lambda_1$.

Knowing the existence of this spectral gap in the normal space, one can 
show that these spectral properties
imply the local stability of the curve $M$ for the nonlinear dynamics given by \eqref{edp:rotators}
(see \cite{cf:Henry} for a general proof, or \cite{cf:GPP} for a proof in our particular model).
In other words, if $p_0$ is sufficiently close to $M$, then there exists a phase $\psi$ such
that $p_t\rightarrow_{t\rightarrow \infty} q_\psi$. This shows that the synchronization is a 
stable phenomenon when $K>1$. But this does not mean that $q_\psi$ alone is stable: 
the trajectory
starting from some perturbation $q_\psi+\epsilon u_t$ of the profile $q_\psi$ may
converge to a profile $q_{\tilde \psi}$ with a phase $\tilde \psi \neq \psi$ (but converging
to $\psi$ when $\epsilon$ goes to $0$).

\medskip

One can in fact do far better than proving the local stability of $M$, by 
describing completely the dynamics.
Using in particular the free energy of the system \cite{cf:BGPoq,cf:GPP}, 
one can prove that, when $K>1$
the solutions of \eqref{edp:rotators} starting from $U=\{p_0\in \mathcal{M}_1(\mathbb{T})
:\, \int_{\mathbb{T}} e^{i\theta}\dd u(\theta)=0\}$
converge to $\frac{1}{2\pi}$, while if $p_0\notin U$ then $p_t$ converges to a $q_\psi\in M$. 

\subsection{Long time behavior}

The end of the preceeding section is devoted to the long time behavior of 
the limit PDE \eqref{edp:rotators} of the model, that is the limit in time of 
the model when the size of the population has already been sent to infinity.
It is the long time behavior of a deterministic system, given by the limit PDE. 
But when $N$ is large but finite, there is still noise in the system. So, in that case, even if
the model is very close to the limit PDE for some finite interval
of time, its behavior can differ dramatically from the one
of the limit PDE in very long times, due to the presence of this noise.

\medskip

Since the noise present in the system disappears when $N$ goes to infinity, one can
see in a certain sense the model with $N$ large as a noisy perturbation of the 
limit PDE. So when $K>1$ and
the empirical measure $\mu_{N,t}$ is close to a synchronized profile $q_\psi\in M$
the dynamics of the finite sized model is given by a competition
between the noise and the operator $L_{q_\psi}$, since this
operator dominates the limit dynamics in the neighborhood of $q_\psi$.
$L_{q_\psi}$ induces a negative
feedback on the normal space,
so in this direction
the process has the same behavior as an
infinite version of an Ornstein Uhlenbeck process
$\dd p_t=-\lambda p_t + \sqrt{\sigma} \dd B_t$
with $\sigma$ small.
So the empirical measure can not go far away in this direction,
and stays close
to $M$.
On the other hand
$L_{q_{\psi}}$ has no effect on the tangent space of $M$
at $q_\psi$, and thus the process can diffuse in this direction.
This means that when $N$ is large $\mu_{N,t}$ stays close to some $q_{\psi^N_t}$, $\psi^N_t$
being some random phase trajectory.
Of course this random phase $\psi^N_t$ converges
to a constant on finite time intervals (since the process converges to a solution of the limit PDE)
and thus to see a macroscopic effect of the noise in the limit one needs to rescale the time.
Since the noise is of size $\sqrt{t/N}$, the appropriate renormalization
is to look at times of order $N$.
It is the purpose of the following Theorem, which corresponds to
Theorem 1.1 in \cite{cf:BGPoq}. We denote
$\Vert\cdot\Vert_{-1}$ the $H^{-1}$ norm.

\begin{thrm}\label{th:long time}
Suppose that $K>1$ and let $\tau_f>0$ and $\psi_0\in \mathbb{T}$.
If for all $\epsilon>0$
\begin{equation}
\lim_{N \to \infty}
\mathbb{P} \left(
\left \Vert \mu_{N, 0} - q_{\psi_0} \right \Vert_{-1}\le \epsilon \right) \, =\, 1\, ,
\end{equation}
then there exists a continuous process $(\psi^N_{\tau})_{\tau \ge 0}$
adapted to the filtration generated by the sequence $(B^j_{N\cdot })_{j=1,\ldots,N}$
such that for all $\epsilon>0$
\begin{equation}
\lim_{N \to \infty}
\mathbb{P} \left( 
\sup_{\tau \in [0, \tau_f]} 
\left\Vert \mu_{N, \tau N} - q_{\psi^N_{\tau} } \right \Vert_{-1}\le \epsilon \right) \, =\, 1\, ,
\end{equation}
and such that $\psi^N_t$ converges weakly to $\psi_0+D_K W_t$, where $(W_t)_{t\ge 0}$ 
is a standard Brownian motion and $D_K$ is a constant that can be computed explicitly 
(in terms of the positive solution $r_K$ of the fixed point problem \eqref{fixed point}).
\end{thrm}

The proof given in \cite{cf:BGPoq} of this Theorem is based on a discretization of 
the dynamics on an intermediate
time-scale, and on a projection of the empirical measure on $M$ at each time step
to follow the fluctuations of the center of synchronization induced by the noise.
This procedure is inspired from the works \cite{cf:BDMP,cf:BBB}, where the authors
show the diffusive behavior of the phase boundary for a one dimensional
reaction-diffusion model with bistable potential and perturbed
with a white noise.
This discretization done, one of the main difficulties in the proof is then to show
that the dynamics of this discretized phase dynamics does not
contain any drift term in the limit $N\rightarrow\infty$, and this is obtained
with use of the symmetry properties of the model: its invariance by rotation has
already been pointed out before, but one can also remark that if 
$p_t(\theta)$ is a solution of \eqref{edp:rotators}, then it is also
the case of $p_t(-\theta)$. See \cite{cf:BGPoq} for more details.

\section{Wasserstein stability of traveling waves for scalar nonlinear advection-diffusion equations}

This section addresses the long time behaviour of the scalar nonlinear advection-diffusion equation
\begin{equation}\label{eq:advd}
  \partial_t u + \partial_x(B(u)) = \frac{1}{2}\partial_x^2(A(u)), \qquad t \geq 0, \quad x \in \xR,
\end{equation}
where $A$ and $B$ are $\xCone$ functions, and $A'(u) = \sigma^2(u) \geq 0$. 
In particular, we aim at illustrating, and extending to Wasserstein distances, the 
classical $\xLone$ stability results of traveling waves of~\eqref{eq:advd}, which 
go back to Osher and Ralston~\cite{OsheRals} as well as Freist\"uhler and 
Serre~\cite{FreiSerr}. We will use probabilistic arguments, based on the interpretation 
by Jourdain and coauthors~\cite{Jour00, JourMalr, JourReyg} of~\eqref{eq:advd} as 
the Fokker-Planck equation of a nonlinear diffusion process. The connection between 
the long time behaviour of this process and the traveling waves of~\eqref{eq:advd} 
was recently pointed out to the author of this survey by Julien Vovelle, to whom warm 
thanks are due.

\subsection{Traveling waves and stationary solutions to~\eqref{eq:advd}} 

When $\sigma^2$ is constant, the equation~\eqref{eq:advd} is a viscous scalar 
conservation law. A (bounded) traveling wave for this equation is a function $\phi$ solving
\begin{equation}\label{eq:phi}
  \frac{\sigma^2}{2}\phi' = B(\phi)-s\phi-q, \qquad \lim_{x \to \pm \infty} \phi(x) = w^{\pm},
\end{equation}
where the speed of the wave $s$ and $w^- \not= w^+$ satisfy the Rankine-Hugoniot condition
\begin{equation}\label{eq:RanHug}
  s = \frac{B(w^+)-B(w^-)}{w^+-w^-},
\end{equation}
which implies $q = B(w^{\pm})-sw^{\pm}$. In the sequel, the boundedness will 
be implicitly assumed and therefore $\phi$ will only be referred to as a {\em traveling wave}.

By the Cauchy-Lipschitz Theorem, for a solution to~\eqref{eq:phi} to exist, it is 
necessary and sufficient that $q \not= B(w)-sw$ for $w$ strictly between $w^-$ 
and $w^+$, which is usually called the Ole{\u\i}nik {\em E}-condition~\cite{Olei}. 
Under this condition, all the traveling waves are translations of each other. If 
$w^- < w^+$, the Ole{\u\i}nik {\em E}-condition rewrites
\begin{equation}
  \forall w \in (w^-, w^+), \qquad \frac{B(w)-B(w^-)}{w-w^-} > \frac{B(w^+)-B(w^-)}{w^+-w^-},
\end{equation}
that is to say the graph of $B$ remains strictly above the line segment joining 
the points $(w^-,B(w^-))$ and $(w^+,B(w^+))$, and it implies that $\phi' > 0$ on 
$\xR$, so that $\phi$ increases from $w^-$ to $w^+$. If $w^->w^+$, then the graph 
of $B$ must remain strictly below this line segment, and $\phi$ decreases from $w^-$ to $w^+$.

If $\phi$ is a traveling wave of~\eqref{eq:advd}, it is immediate that $u(t,x) := \phi(x-st)$ 
solves~\eqref{eq:advd}. If $s=0$, then $u$ is a stationary solution to~\eqref{eq:advd}. 
It has been known since the works by Il'in and Ole{\u\i}nik~\cite{IlinOlei58, IlinOlei60} 
that traveling waves describe the long time behavior of solutions to~\eqref{eq:advd}. 
The following $\xLone$ stability theorem is due Freist\"uhler and Serre~\cite{FreiSerr} 
and is based on a former result by Osher and Ralston~\cite{OsheRals}.

\begin{thrm}[Freist\"uhler and Serre]\label{thrm:FreiSerr}
  Let $\phi$ be a traveling wave of~\eqref{eq:advd} and $u_0$ such that 
  $u_0 - \phi \in \xLone(\xR)$. The solution $u(t,x)$ of~\eqref{eq:advd} with initial 
  datum $u_0$ satisfies
  \begin{equation}
    \lim_{t \to +\infty} ||u(t,\cdot) - \phi(\cdot - st + \delta)||_{\xLone(\xR)} = 0,
  \end{equation}
  where the phase shift $\delta$ is defined by
  \begin{equation}\label{eq:delta}
    \delta := \frac{1}{w^+-w^-} \int_{x \in \xR} (u_0(x)-\phi(x)) \xdif x.
  \end{equation}
\end{thrm}
Note that the definition~\eqref{eq:delta} of $\delta$ ensures that
\begin{equation}\label{eq:delta2}
  \int_{x \in \xR} (u_0(x)-\phi(x+\delta)) \xdif x = 0.
\end{equation}

For a nonconstant diffusion coefficient $\sigma^2$, the equation~\eqref{eq:phi} 
defining a traveling wave has to be replaced with
\begin{equation}\label{eq:phiNC}
  \frac{1}{2}(A(\phi))' = B(\phi)-s\phi-q, \qquad \lim_{x \to \pm \infty} \phi(x) = w^{\pm},
\end{equation}
where $s$ is still defined in terms of $w^- \not= w^+$ by the Rankine-Hugoniot 
condition~\eqref{eq:RanHug}. Under the uniform ellipticity condition $\inf_u \sigma^2(u) > 0$, 
the Ole{\u\i}nik {\em E}-condition remains necessary and sufficient for a traveling wave 
to exist, and traveling waves remain monotonic on the real line. Gasnikov~\cite{Gasn} 
proved that, if $A$ and $B$ are $\xCn{4}$ on $[w^-,w^+]$ (or $[w^+,w^-]$), then 
Theorem~\ref{thrm:FreiSerr} holds without any change in its statement.

\subsection{Probabilistic interpretation of~\eqref{eq:advd}}\label{ss:probinterp}

Traveling waves to~\eqref{eq:advd} with $w^-=0$ and $w^+=1$ can be interpreted 
as cumulative distribution functions (CDFs) of probability measures on the real line. 
If we assume that the initial datum $u_0$ of~\eqref{eq:advd} is also the CDF of a 
probability measure $m$, then $u(t,\cdot)$ remains the CDF of a probability measure 
$P_t$ at all times. Besides, taking the space derivative of~\eqref{eq:advd} yields 
the formal evolution equation
\begin{equation}\label{eq:efp}
  \partial_t P_t = \frac{1}{2}\partial_x^2\left(\sigma^2(H*P_t(x)) P_t\right) 
  - \partial_x \left(b(H*P_t(x)) P_t\right)
\end{equation}
for $P_t$, where $b := B'$ and $H*P_t = u(t, \cdot)$ refers to the spatial convolution 
of $P_t$ with the Heaviside function. The equation~\eqref{eq:efp} is the Fokker-Planck 
equation of the diffusion process
\begin{equation}
  \left\{\begin{aligned}
    & \xdif X_t = b(H*P_t(X_t)) \xdif t + \sigma(H*P_t(X_t)) \xdif W_t,\\
    & \text{$P_t$ is the law of $X_t$},
  \end{aligned}\right.\end{equation}
where $(W_t)_{t \geq 0}$ is a standard real Brownian motion, and $X_0$ is distributed 
according to $m$, independently of $(W_t)_{t \geq 0}$. Note that the coefficients of 
this stochastic differential equation depend on the law $P_t$ of $X_t$, which is the trace 
of the nonlinearity of the Fokker-Planck equation~\eqref{eq:efp}. Therefore, the process 
$(X_t)_{t \geq 0}$ is said to be {\em nonlinear in McKean's sense}.

The existence and uniqueness of the nonlinear process $(X_t)_{t \geq 0}$ were 
established in~\cite{JourReyg} under the assumptions that $b$ and $\sigma^2$ 
are continuous on $[0,1]$, $m$ have a finite first order moment, 
$\sigma^2(u) > 0$ on $(0,1)$ and:
\begin{itemize}
  \item if $\sigma^2(0)=0$, then $u_0(x) > 0$ for all $x \in \xR$;
  \item if $\sigma^2(1)=0$, then $u_0(x) < 1$ for all $x \in \xR$;
\end{itemize}
where $u_0 := H*m$ is the CDF of the initial distribution $m$. 

\begin{rmrk}\label{rmrk:CDF}
  The assumption on the finiteness of the first order moment of $m$ is natural in 
  order to obtain $\xLone$ stability results on the solution to~\eqref{eq:advd}. Indeed, 
  in general, if $F_1$ and $F_2$ are the CDFs of probability measures $m_1$ and 
  $m_2$ on $\xR$, then $||F_1-F_2||_{\xLone(\xR)}$ need not be finite, but if we 
  assume in addition that the first order moments of $m_1$ and $m_2$ are finite, 
  then $||F_1-F_2||_{\xLone(\xR)} < +\infty$, and the difference between the 
  expectations of $m_1$ and $m_2$ is given by
  \begin{equation}
    \int_{x \in\R} x m_1(\xdif x) - \int_{x \in\R} x m_2(\xdif x) = \int_{x \in \R} (F_1(x)-F_2(x))\xdif x.
  \end{equation}
\end{rmrk}

We first provide a probabilistic interpretation of the speed $s$ of a traveling wave 
as the average velocity of the nonlinear process $(X_t)_{t \geq 0}$. Indeed, the 
expectation of $X_t$ satisfies
\begin{equation}
  \Exp[X_t] = \Exp[X_0] + \int_{s=0}^t \Exp[b(H*P_s(X_s))]\xdif s.
\end{equation}
Besides, it was proved in~\cite{JourReyg} that, $\xdif s$-almost everywhere, the 
measure $P_s$ does not weight points, which implies that $H*P_s(X_s)$ is 
uniformly distributed on $[0,1]$. We therefore rewrite
\begin{equation}\label{eq:ExpXt}
  \Exp[X_t] = \Exp[X_0] + \int_{s=0}^t \int_{u=0}^1 b(u) \xdif u \xdif s = \Exp[X_0] + st,
\end{equation}
where $s$ is given by the Rankine-Hugoniot condition $s = B(1)-B(0)$. 

We now describe the long time behaviour of the nonlinear process in terms of 
traveling waves. We first discuss conditions ensuring that traveling waves are 
well defined thanks to the following lemma, which was obtained 
in~\cite[Proposition~4.1 and Corollary~4.4]{JourReyg} by solving~\eqref{eq:phiNC} explicitly.

\begin{lmm}\label{lmm:phi}
  Assume that $b$ and $\sigma^2$ are continuous on $[0,1]$, that $A$ is increasing 
  on $[0,1]$, and that the Ole{\u\i}nik {\em E}-condition
  \begin{equation}\label{eq:EB}
    \forall u \in (0,1), \qquad B(u) > B(0) + su
  \end{equation}
  is satisfied, where $s$ is defined by the Rankine-Hugoniot condition $s = B(1)-B(0)$. 
  
  Then there exists a traveling wave $\phi$ increasing from $0$ to $1$, and all such 
  traveling waves are translations of each other. Besides, the probability measure 
  with CDF $\phi$ has a finite first order moment if and only if
  \begin{equation}\label{eq:mom}
    \int_{u=0}^{1/2} \frac{u\sigma^2(u)}{B(u)-B(0)-su} 
    \xdif u + \int_{u=1/2}^1 \frac{(1-u)\sigma^2(u)}{B(u)-B(0)-su} \xdif u < +\infty.
  \end{equation}
\end{lmm}
Let us mention that, if the Ole{\u\i}nik {\em E}-condition is relaxed by allowing that 
$B(u)=B(0)+su$ for some $u \in (0,1)$ such that $\sigma^2(u)=0$, then one can 
exhibit traveling waves that are not translations of each other, 
see~\cite[Remark~4.2]{JourReyg}.

Let us now fix a probability measure $m$ with finite first order moment, and let $\phi$ 
be given by Lemma~\ref{lmm:phi}. If~\eqref{eq:mom} holds, then by Remark~\ref{rmrk:CDF}, 
$u_0-\phi \in \xLone(\xR)$; besides, choosing the phase shift $\delta$ in order to 
satisfy~\eqref{eq:delta2} amounts to selecting $u_{\infty} = \phi(\cdot + \delta)$ having 
the same expectation as $u_0$. By~\eqref{eq:ExpXt}, we already know that the expectation 
of $X_t-st$ is constant and equal to the expectation of $u_{\infty}$. 
Theorem~\ref{thrm:FreiSerr} contains the much stronger statement that the long 
time behaviour of this process is described by the stationary wave $u_{\infty}$, in the 
sense that the CDF of $X_t-st$ converges to $u_{\infty}$ in $\xLone(\xR)$. In the next 
subsection, we extend this result to Wasserstein distances.

\subsection{Contraction and convergence to equilibrium in Wasserstein distance}

We recall that the Wasserstein distance of order $p \in [1,+\infty)$ between 
two probability measures on the real line with respective CDFs $F_1$ and $F_2$ is given by
\begin{equation}
  \Wass_p(F_1, F_2) := \left(\int_{w=0}^1 |F_1^{-1}(w) - F_2^{-1}(w)|^p \xdif w\right)^{1/p},
\end{equation}
where the pseudo-inverse $F^{-1}$ of a CDF $F$ is defined by 
$F^{-1}(w) := \inf\{x \in \xR: F(x) \geq w\}$ for all $w \in (0,1)$. Note that, in particular,
\begin{equation}\label{eq:W1L1}
  \Wass_1(F_1, F_2) = ||F_1 - F_2||_{\xLone(\xR)}.
\end{equation}

Our first result is the following Wasserstein contraction property of~\eqref{eq:advd}. 
Given two CDFs $u_0$, $v_0$ on the real line, we now denote by $u_t := u(t,\cdot)$ 
and $v_t := v(t,\cdot)$ the corresponding solutions to~\eqref{eq:advd}.
\begin{prpstn}\label{prpstn:contract}
  Assume that $b$ and $\sigma^2$ are continuous on $[0,1]$, that $A$ is increasing 
  on $[0,1]$ and that $m$ has a finite first order moment. For all $p \in [1,+\infty)$,
  \begin{itemize}
    \item if $\Wass_p(u_0,v_0) = +\infty$, then $\Wass_p(u_t, v_t) = +\infty$ for all $t \geq 0$;
    \item if $\Wass_p(u_0,v_0) < +\infty$, then $t \mapsto \Wass_p(u_t, v_t)$ is 
    nonincreasing on $[0,+\infty)$.
  \end{itemize}
\end{prpstn}

The proof of Proposition~\ref{prpstn:contract} is detailed in~\cite[Proposition~3.1]{JourReyg}. 
It is entirely probabilistic, and relies on a coupling argument for the order statistics of a 
system of mean-field interacting particles approximating $u_t$ and $v_t$. We note that, 
by~\eqref{eq:W1L1}, the case $p=1$ of Proposition~\ref{prpstn:contract} is nothing but 
the classical $\xLone$ stability estimate
\begin{equation}
  \forall t \geq 0, \qquad ||u_t - v_t||_{\xLone(\xR)} \leq ||u_0 - v_0||_{\xLone(\xR)}
\end{equation}
for~\eqref{eq:advd}. Similar Wasserstein estimates were obtained for scalar conservation 
laws, that is to say $\sigma^2=0$ in our setting, by Bolley, Brenier and Loeper in~\cite{bolley}.

Assuming classical regularity for $u$ and $v$, one can actually go deeper into the 
description of the evolution of $\Wass_p(u_t, v_t)$. Indeed, observing that the 
pseudo-inverse $u_t^{-1}(w)$ of the CDF $u_t$ satisfies the equation
\begin{equation}
  \partial_t u_t^{-1} = b(w) - \partial_w \left(\frac{\sigma^2(w)}{2\partial_w u_t^{-1}}\right),
\end{equation}
one can derive the explicit formula for the time derivative of $\Wass_p(u_t, v_t)$:
\begin{equation}
  \xDrv{\Wass_p(u_t, v_t)^p}{t} = 
  -\frac{p(p-1)}{2} \int_{w=0}^1 \sigma^2(w) |u_t^{-1} - v_t^{-1}|^{p-2} 
  \frac{\left(\partial_w u_t^{-1} - \partial_w v_t^{-1}\right)^2}{\partial_w u_t^{-1}
  \partial_w v_t^{-1}}\xdif w.
\end{equation}
This leads to the following convergence theorem, which is the main result 
of~\cite{JourReyg} and readily extends Theorem~\ref{thrm:FreiSerr} to Wasserstein distances.
\begin{thrm}\label{thrm:JourReyg}
  Assume that:
  \begin{itemize}
    \item $b$ is $\xCn{{1+\beta}}$ on $[0,1]$, $\sigma^2$ is positive and 
    $\xCn{{2+\beta}}$ on $[0,1]$, and the equilibrium conditions~\eqref{eq:EB} 
    and~\eqref{eq:mom} hold;
    \item the probability measure $m$ has a finite first order moment;
    \item the Wasserstein distance of order $2$ between $u_0 := H*m$ and any 
    traveling wave $\phi$ increasing from $0$ to $1$ is finite.
  \end{itemize}
  Let us denote by $u_{\infty}$ the traveling wave with the same expectation 
  as $u_0$. Then, for all $p \geq 2$ such that $\Wass_p(u_0, u_{\infty}) < +\infty$, we have
  \begin{equation}
    \forall 1 \leq q < p, \qquad \lim_{t \to +\infty} \Wass_q(u(t, \cdot), u_{\infty}(\cdot - st)) = 0.
  \end{equation}
\end{thrm}

\subsection{Conclusion}

We have interpreted the scalar nonlinear advection-diffusion equation~\eqref{eq:advd}, 
with a CDF as an initial datum, as the Fokker-Planck equation of a nonlinear diffusion 
process $(X_t)_{t \geq 0}$. The expectation of this process evolves linearly in time, at 
a velocity given by the speed of traveling waves of~\eqref{eq:advd} increasing from $0$ 
to $1$. Under the Ole{\u\i}nik {\em E}-condition~\eqref{eq:EB}, Theorem~\ref{thrm:FreiSerr} 
shows that the fluctuation of $X_t$ around $st$ converges, in $\xLone(\xR)$, to an 
equilibrium distribution described by the traveling wave having the same expectation 
as $X_t-st$. Theorem~\ref{thrm:JourReyg} extends this result to Wasserstein distance.

The probabilistic interpretation of~\eqref{eq:advd} can also lead to further developments 
on Theorems~\ref{thrm:FreiSerr} and~\ref{thrm:JourReyg}. For example, in the case of 
a constant diffusion coefficient $\sigma^2$, an exponential rate of decay to equilibrium 
for $X_t-st$ was obtained in~\cite{JourMalr}, for initial distributions close to the equilibrium 
distribution. The decay was expressed in $\chi_2$ distance, and using the transport 
chi-square inequality of~\cite{Jour12}, it can be translated in quadratic Wasserstein 
distance. We refer to~\cite[Subsection~4.3]{JourReyg} for details in this direction.

\bibliographystyle{plain}
\bibliography{session-Markov}

\begin{thebibliography}{10}

\bibitem{BCGMZ13}
J.-B. Bardet, A.~Christen, A.~Guillin, F.~Malrieu, and P.-A. Zitt.
\newblock Total variation estimates for the {TCP} process.
\newblock {\em Electron. J. Probab}, 18(10):1--21, 2013.

\bibitem{BCT08}
P.~Bertail, S.~Cl{\'e}men{\c{c}}on, and J.~Tressou.
\newblock A storage model with random release rate for modeling exposure to
  food contaminants.
\newblock {\em Math. Biosci. Eng.}, 5(1):35--60, 2008.

\bibitem{cf:BBB}
L.~Bertini, S.~Brassesco, and P.~Butt\`a.
\newblock Soft and hard wall in a stochastic reaction diffusion equation.
\newblock {\em Arch. Ration. Mech. Anal.}, 190:307--345, 2008.

\bibitem{cf:BGP}
L.~Bertini, G.~Giacomin, and K.~Pakdaman.
\newblock Dynamical aspects of mean field plane rotators and the {K}uramoto
  model.
\newblock {\em Journal of Statistical Physics}, 138:270--290, 2010.

\bibitem{cf:BGPoq}
L.~Bertini, G.~Giacomin, and C.~Poquet.
\newblock Synchronisation and random long time dynamics for mean-field plane
  rotators.
\newblock {\em Probab. Theory Relat. Fields}, 2013.

\bibitem{bolley}
F.~Bolley, Y.~Brenier, and G.~Loeper.
\newblock Contractive metrics for scalar conservation laws.
\newblock {\em J. Hyperbolic Differ. Equ.}, 2(1):91--107, 2005.

\bibitem{Bou14}
F.~Bouguet.
\newblock Quantitative speeds of convergence for exposure to food contaminants.
\newblock Preprint available on ar{X}iv, to appear in ESAIM P\& S, 2013.

\bibitem{cf:BDMP}
S.~Brassesco, A.~De~Masi, and E.~Presutti.
\newblock Brownian fluctuations of the interface in the $d=1$
  {G}inzburg-{L}andau equation with noise.
\newblock {\em Annal. Inst. H. Poincar\'e}, 31:81--118, 1995.

\bibitem{Dav93}
M.~H.~A. Davis.
\newblock {\em Markov models and optimization}, volume~49 of {\em Monographs on
  Statistics and Applied Probability}.
\newblock Chapman \& Hall, London, 1993.

\bibitem{DownMeynTweedie}
D.~Down, S.~P. Meyn, and R.~L. Tweedie.
\newblock Exponential and uniform ergodicity of {M}arkov processes.
\newblock {\em Ann. Probab.}, 23(4):1671--1691, 1995.

\bibitem{fontbona-panloup}
J.~Fontbona and F.~Panloup.
\newblock Rate of convergence to equilibrium of fractional driven stochastic
  differential equations with some multiplicative noise.
\newblock {\em Preprint available at
  http://hal.archives-ouvertes.fr/hal-00989414}, pages 1--34, 2014.

\bibitem{FreiSerr}
H.~Freist{{\"u}}hler and D.~Serre.
\newblock {$L^1$} stability of shock waves in scalar viscous conservation laws.
\newblock {\em Comm. Pure Appl. Math.}, 51(3):291--301, 1998.

\bibitem{GPS}
S.~Gadat, F.~Panloup, and S.~Saadane.
\newblock Regret bounds for narendra-shapiro bandit algorithms.
\newblock Preprint available on ar{X}iv n° 1502.04874, 2015.

\bibitem{cf:Gartner}
J.~G\"artner.
\newblock On {M}c{K}ean-{V}lasov limit for interacting diffusions.
\newblock {\em Math. Nachr.}, 137:197--248, 1988.

\bibitem{Gasn}
A.~V. Gasnikov.
\newblock Time-asymptotic behavior of the solution of the initial {C}auchy
  problem for a conservation law with nonlinear divergent viscosity.
\newblock {\em Izv. Ross. Akad. Nauk Ser. Mat.}, 73(6):39--76, 2009.

\bibitem{cf:GPP}
G.~Giacomin, K.~Pakdaman, and X.~Pellegrin.
\newblock Global attractor and asymptotic dynamics in the {K}uramoto model for
  coupled noisy phase oscillators.
\newblock {\em Nonlinearity}, 25:1247--1273, 2012.

\bibitem{hairer}
M.~Hairer.
\newblock Ergodicity of stochastic differential equations driven by fractional
  {B}rownian motion.
\newblock {\em Ann. Probab.}, 33(2):703--758, 2005.

\bibitem{hairer2}
M.~Hairer and A.~Ohashi.
\newblock Ergodic theory for {SDE}s with extrinsic memory.
\newblock {\em Ann. Probab.}, 35(5):1950--1977, 2007.

\bibitem{hairer-pillai}
M.~Hairer and N.~S. Pillai.
\newblock Regularity of laws and ergodicity of hypoelliptic {SDE}s driven by
  rough paths.
\newblock {\em Ann. Probab.}, 41(4):2544--2598, 2013.

\bibitem{cf:Henry}
D.~Henry.
\newblock {\em Geometric theory of semilinear parabolic equations}, volume 840
  of {\em Lecture Notes in Mathematics}.
\newblock Springer-Verlag, 1981.

\bibitem{IlinOlei58}
A.~M. Il'in and O.~A. Ole{\u\i}nik.
\newblock Behavior of solutions of the {C}auchy problem for certain quasilinear
  equations for unbounded increase of the time.
\newblock {\em Dokl. Akad. Nauk SSSR}, 120:25--28, 1958.

\bibitem{IlinOlei60}
A.~M. Il'in and O.~A. Ole{\u\i}nik.
\newblock Asymptotic behavior of solutions of the {C}auchy problem for some
  quasi-linear equations for large values of the time.
\newblock {\em Mat. Sb. (N.S.)}, 51 (93):191--216, 1960.

\bibitem{Jour00}
B.~Jourdain.
\newblock Diffusion processes associated with nonlinear evolution equations for
  signed measures.
\newblock {\em Methodol. Comput. Appl. Probab.}, 2(1):69--91, 2000.

\bibitem{Jour12}
B.~Jourdain.
\newblock Equivalence of the {P}oincar\'e inequality with a
  transport-chi-square inequality in dimension one.
\newblock {\em Electron. Commun. Probab.}, 17:no. 43, 12, 2012.

\bibitem{JourMalr}
B.~Jourdain and F.~Malrieu.
\newblock Propagation of chaos and {P}oincar{\'e} inequalities for a system of
  particles interacting through their {CDF}.
\newblock {\em Ann. Appl. Probab.}, 18(5):1706--1736, 2008.

\bibitem{JourReyg}
B.~Jourdain and J.~Reygner.
\newblock Propagation of chaos for rank-based interacting diffusions and long
  time behaviour of a scalar quasilinear parabolic equation.
\newblock {\em Stochastic Partial Differential Equations: Analysis and
  Computations}, 1(3):455--506, 2013.

\bibitem{LP08b}
D.~Lamberton and G.~Pag{\`e}s.
\newblock How fast is the bandit?
\newblock {\em Stoch. Anal. Appl.}, 26(3):603--623, 2008.

\bibitem{LP08}
D.~Lamberton and G.~Pag{\`e}s.
\newblock A penalized bandit algorithm.
\newblock {\em Electron. J. Probab.}, 13:no. 13, 341--373, 2008.

\bibitem{LPT04}
D.~Lamberton, G.~Pag{\`e}s, and P.~Tarr{\`e}s.
\newblock When can the two-armed bandit algorithm be trusted?
\newblock {\em Ann. Appl. Probab.}, 14(3):1424--1454, 2004.

\bibitem{Lin02}
T.~Lindvall.
\newblock {\em Lectures on the coupling method}.
\newblock Dover Publications Inc., Mineola, NY, 2002.
\newblock Corrected reprint of the 1992 original.

\bibitem{Nualart02}
D.~Nualart and A.~R{\u{a}}{\c{s}}canu.
\newblock Differential equations driven by fractional {B}rownian motion.
\newblock {\em Collect. Math.}, 53(1):55--81, 2002.

\bibitem{Olei}
O.~A. Ole{\u\i}nik.
\newblock Uniqueness and stability of the generalized solution of the {C}auchy
  problem for a quasi-linear equation.
\newblock {\em Uspehi Mat. Nauk}, 14(2 (86)):165--170, 1959.

\bibitem{OsheRals}
S.~Osher and J.~Ralston.
\newblock {$L^{1}$} stability of travelling waves with applications to
  convective porous media flow.
\newblock {\em Comm. Pure Appl. Math.}, 35(6):737--749, 1982.

\bibitem{cf:pikovsky}
A.S. Pikovsky and J.~Kurths.
\newblock Coherence resonance in a noise driven excitable system.
\newblock {\em Phys. Rev. Lett.}, 78:775--778, 1997.

\bibitem{ruadulescu}
S.~R{\u{a}}dulescu and M.~R{\u{a}}dulescu.
\newblock An application of {H}adamard-{L}\'evy's theorem to a scalar initial
  value problem.
\newblock {\em Proc. Amer. Math. Soc.}, 106(1):139--143, 1989.

\bibitem{cf:Sznitman}
A.-S. Sznitman.
\newblock Topics in propagation of chaos.
\newblock In {\em \'Ecole d'\'et\'e de probabilit\'es de {S}aint-{F}lour
  XIX--1989}, volume 1464 of {\em Lecture Notes in Math.} Springer, 1991.

\end{thebibliography}
\end{document}